\documentclass[12pt,letterpaper]{amsart}
\usepackage[english]{babel}
\usepackage[dvips]{graphicx}
\usepackage{amscd}
\usepackage{latexsym}
\usepackage{amsfonts}
\newtheorem{theorem}{Theorem}
\newtheorem{proposition}[theorem]{Proposition}
\newtheorem{lemma}[theorem]{Lemma}
\newtheorem{corollary}[theorem]{Corollary}
\begin{document}
\title{The Strong Symmetric Genus of the Finite Coxeter Groups}
\author{Michael A. Jackson}
\begin{abstract}
The strong symmetric genus of a finite group $G$ is the smallest genus of a closed orientable 
topological surface on which $G$ acts faithfully as a group of orientation preserving automorphisms. In this paper we complete 
the calculation of the strong symmetric genus for each finite Coxeter group excluding the group $E_8$.
\end{abstract}
\maketitle
\section{Introduction}

It is known that the automorphism group of a compact Riemann surface of 
genus $g \ge 2$ is finite; furthermore, the size of this automorphism group must be no larger 
than $84(g-1)$ when $ g\ge 2$ \cite{h:uag}. When a finite group $G$ can be realized as 
an automorphism group of a genus $g$ compact Riemann surface such that $|G|=84(g-1)$, 
we say that $G$ is a Hurwitz group. 
If a group $G$ is not the automorphism group of any Riemann surface of genus 0 or 1, then for any 
Riemann surface on which the group acts as an automorphism group the genus must be at least 
$1+\frac{|G|}{84}$. Burnside \cite{b:tgfo} began the investigation of a related problem:  
find the the least genus $g$ of a Riemann surface on which a given finite group acts faithfully as a group 
of automorphisms. Equivalently, given a finite group $G$, one may want to find the least genus $g$ 
of a closed orientable topological surface 
on which $G$ acts as a group of orientation preserving symmetries.  
The latter parameter is denoted $\sigma^0(G)$ and has become known as the strong symmetric
genus of $G$ (see \cite[Chapter 6]{gt:tgt} ). 

Many results are known concerning the 
strong symmetric genus of finite groups. All finite groups with strong symmetric genus 
less than 4 are known \cite{b:cfga,mz:gsssg}. It has been shown as well that for each positive integer $n$, there is a finite group 
$G$ with $\sigma ^{0} (G)=n$ \cite{mz:essg} . 
The strong symmetric genus of several infinite families of finite groups 
have been found: the alternating and symmetric groups \cite{c:gasg,c:mgasg,c:srqtg}, 
the hyperoctahedral groups \cite{j:ssgh}, the 
groups $PSL_2 (q)$ \cite{gs:gop,gs:rpsl}, and the groups $SL_2(q)$ \cite{v:gos}. In addition, the 
strong symmetric genus has been found for the sporadic finite simple groups 
\cite{cww:sgsg,w:mhg,w:sgfg,w:sgbm}. 

The symmetric groups and the hyperoctahedral groups are two infinite families of finite 
Coxeter groups. They are often referred to as the $A$-type and $B$-type Coxeter groups, respectively. 
As stated above, the strong symmetric genus is known for each group in these families. Another family of 
finite Coxeter groups is the dihedral groups, which are automorphism groups of the sphere and, thus, 
have strong symmetric genus 0. This leaves one infinite family of finite Coxeter groups whose 
strong symmetric genus has not been found previously: 
the $D$-type groups. In this paper we calculate the strong symmetric genus 
of the $D$-type finite Coxeter groups and of the sporadic finite Coxeter groups, excluding $E_8$. These 
new results will be given in Theorem \ref{thm:main} below; the previously known 
results concerning the other finite Coxeter groups will also be stated in Theorem \ref{thm:main}.

The strong symmetric genus of the finite Coxeter groups will be shown 
by demonstrating the existence of certain pairs of generators. 
If a finite group $G$ has generators $x$ and $y$ of orders $p$ and $q$, respectively, with $xy$ having 
order $r$, then we say that $(x,y)$ is a $(p,q,r)$ generating pair of $G$. By the obvious symmetries concerning 
generators, we will use the convention that $p\leq q\leq r$. Following the convention of Marston Conder 
\cite{c:srqtg}, we say that a $(p,q,r)$ generating pair of $G$ is a minimal generating pair if there does not 
exist a $(k,l,m)$ generating pair for $G$ with $\frac{1}{k}+\frac{1}{l}+\frac{1}{m} >\frac{1}{p} +\frac{1}{q} +
\frac{1}{r}$. 

Before stating Theorem \ref{thm:main}, we define the groups involved. For an $n\geq 3$, $B_n$ will be the 
group of symmetries of the $n$ dimensional cube, and $D_n$ will be the group of orientation preserving symmetries
of the $n$ dimensional cube. For the sporadic finite Coxeter groups, we list the results in Table \ref{tab:spdef}.
\begin{table}[h] 
\begin{tabular}{|c|c|c|c|}\hline 
group $G$ & description \\
\hline
\hline
$G_2$ & isomorphic to the dihedral group of order 12 \\
\hline
$I_3$ & symmetry group of the regular icosahedron \\
\hline
$I_4$ &  symmetry group of the 4- dimensional, regular 120-cell\\
\hline
$F_4$ &  symmetry group of the 4-dimensional, regular 24-cell\\
\hline
$E_6$ &  symmetry group of the $E_6$ polytope\\
\hline
$E_7$ &  symmetry group of the $E_7$ polytope\\
\hline
$E_8$ & symmetry group of the $E_8$ polytope\\
\hline
\end{tabular}
\caption{Sporadic Coxeter group descriptions} \label{tab:spdef}
\end{table}

\begin{theorem}\label{thm:main}
Let $G$ be a finite Coxeter group. If $G$ is the dihedral group of size $2n$, then $G$ has a $(2,2,n)$ 
minimal generating pair and $\sigma ^0 (G)=0$. If $G=\Sigma _n$ for $n>29$, then $G$ has a $(2,3,8)$ 
minimal generating pair and $\sigma ^0 (\Sigma _n) = \frac{n!}{48}$. If $G=B_n$ for $n>8$, then $G$ has a $(2,4,6)$ 
minimal generating pair and $\sigma ^0 (B_n)=\frac{n! \; 2^n}{24}$. If $G=D_n$ for $n>29$, then $G$ has a $(2,3,8)$ 
minimal generating pair and $\sigma ^0 (D_n)=\frac{n! \; 2^{n-1}}{48}$. The results for $G$ being one of the sporadic 
cases are given in Table \ref{tab:spor}. The remaining cases are listed in Table \ref{tab:exc}.
\end{theorem}
\begin{table}[h] 
\begin{tabular}{|c|c|c|c|}\hline 
group $G$ & size & min. gen. pair & $\sigma ^0 (G)$ \\
\hline
\hline
$G_2$ & $12$ & $(2,2,6)$ & $0$ \\
\hline
$I_3$ & $120$ & $(2,3,10)$ & $5$\\
\hline
$I_4$ & $14400$ & $(2,4,6)$ & $601$\\
\hline
$F_4$ & $1152$ & $(2,6,6)$ & $97$\\
\hline
$E_6$ & $51,840$ & $(2,4,8)$ & $3241$\\
\hline
$E_7$ & $2,903,040$ & $(2,4,7)$ & $155,521$\\
\hline
\end{tabular}
\caption{Sporadic Coxeter group results} \label{tab:spor}
\end{table}
\begin{table}[h]
\begin{tabular}{|c|c|c|c|c|c|c|}\hline 
 & \multicolumn{2}{c} {$\Sigma _n$}  \vline & \multicolumn{2}{c}{$B_n$} \vline & \multicolumn{2}{c}{$D_n$} \vline \\
$n$ & best $\Delta $ & $\sigma ^0 (\Sigma _n)$ & best $\Delta $ & $\sigma ^0 (B _n)$ & best $\Delta $ & $\sigma ^0 (D _n)$\\
\hline
\hline 
$n=3$ & \multicolumn{2}{c}{Dihedral of order 6} \vline & $(2,4,6)$ & $\frac{3! 2^3}{24}+1=3$&
 \multicolumn{2}{c}{$D _3 = \Sigma _4$} \vline \\
\hline
$n=4$ & $(2,3,4)$ & $0$  & $(2,4,6)$ & $\frac{4! 2^4}{24}+1=17$& $(3,4,4)$ & $\frac{4! 2^3}{12}+1=17$ \\
\hline
$n=5$ & $(2,4,5)$ & $\frac{5!}{40}+1=4$  & $(2,4,10)$ & $\frac{5! 2^5 (3)}{40}+1$& $(2,4,5)$ & $\frac{5! 2^4}{40}+1$ \\
\hline
$n=6$ & $(2,5,6)$ & $\frac{6!}{15}+1$  & $(2,6,6)$ & $\frac{6! 2^6}{12}+1$& $(2,5,6)$ & $\frac{6! 2^5}{15}+1$ \\
\hline
$n=7$ & $(2,3,10)$ & $\frac{7!}{30}+1$  & $(2,4,6)$ & $\frac{7! 2^7}{24}+1$& $(2,4,6)$ & $\frac{7! 2^6}{24}+1$ \\
\hline
$n=8$ & $(2,4,7)$ & $\frac{8!(3)}{56}+1$  & $(2,4,8)$ & $\frac{8! 2^8}{16}+1$& $(2,4,7)$ & $\frac{8! 2^7 (3)}{56}+1$ \\
\hline
$n=9$ & $(2,4,6)$ & $\frac{9!}{24}+1$  & $(2,4,6)$ & $\frac{9! 2^9}{24}+1$& $(2,4,6)$ & $\frac{9! 2^8}{24}+1$ \\
\hline
$n=10$ & $(2,3,10)$ & $\frac{10!}{30}+1$  & $(2,4,6)$ & $\frac{10! 2^{10}}{24}+1$& $(2,3,10)$ & $\frac{10! 2^9}{30}+1$ \\
\hline
$n=11$ & $(2,4,5)$ & $\frac{11!}{40}+1$  & $(2,4,6)$ & $\frac{11! 2^{11}}{24}+1$& $(2,4,5)$ & $\frac{11! 2^{10}}{40}+1$ \\
\hline
$n=12$ & $(2,3,12)$ & $\frac{12!}{24}+1$  & $(2,4,6)$ & $\frac{12! 2^{12}}{24}+1$& $(2,3,12)$ & $\frac{12! 2^{11}}{24}+1$ \\
\hline
$n=13$ & $(2,3,12)$ & $\frac{13!}{24}+1$  & $(2,4,6)$ & $\frac{13! 2^{13}}{24}+1$& $(2,3,12)$ & $\frac{13! 2^{12}}{24}+1$ \\
\hline
$n=14$ & $(2,4,6)$ & $\frac{14!}{24}+1$  & $(2,4,6)$ & $\frac{14! 2^{14}}{24}+1$& $(2,3,14)$ & $\frac{14! 2^{13}}{21}+1$ \\
\hline
$n=15$ & $(2,4,5)$ & $\frac{15!}{40}+1$  & $(2,4,6)$ & $\frac{15! 2^{15}}{24}+1$& $(2,4,5)$ & $\frac{15! 2^{14}}{40}+1$ \\
\hline
$n=16$ & $(2,4,5)$ & $\frac{16!}{40}+1$  & $(2,4,6)$ & $\frac{16! 2^{16}}{24}+1$& $(2,4,5)$ & $\frac{16! 2^{15}}{40}+1$ \\
\hline
$n=17$ & $(2,4,6)$ & $\frac{17!}{24}+1$  & $(2,4,6)$ & $\frac{17! 2^{17}}{24}+1$& $(2,4,6)$ & $\frac{17! 2^{16}}{24}+1$ \\
\hline
$n=20$ & $(2,3,8)$ & $\frac{20!}{48}+1$  & $(2,4,6)$ & $\frac{20! 2^{20}}{24}+1$& $(2,4,5)$ & $\frac{20! 2^{19}}{40}+1$ \\
\hline
$n=22$ & $(2,3,10)$ & $\frac{22!}{30}+1$  & $(2,4,6)$ & $\frac{22! 2^{22}}{24}+1$& $(2,3,10)$ & $\frac{22! 2^{21}}{30}+1$ \\
\hline
$n=23$ & $(2,3,10)$ & $\frac{23!}{30}+1$  & $(2,4,6)$ & $\frac{23! 2^{23}}{24}+1$& $(2,3,12)$ & $\frac{23! 2^{22}}{24}+1$ \\
\hline
$n=26$ & $(2,4,5)$ & $\frac{26!}{40}+1$  & $(2,4,6)$ & $\frac{26! 2^{26}}{24}+1$& $(2,4,5)$ & $\frac{26! 2^{25}}{40}+1$ \\
\hline
$n=29$ & $(2,3,12)$ & $\frac{29!}{24}+1$  & $(2,4,6)$ & $\frac{29! 2^{29}}{24}+1$& $(2,3,12)$ & $\frac{29! 2^{28}}{24}+1$ \\
\hline
\end{tabular}
\caption{Exceptional case results}  \label{tab:exc}
\end{table}

\section{Generating Pairs and Strong Symmetric Genus} \label{sec:gp}

Recall that the groups of small strong symmetric genus are well known 
(see \cite{b:cfga,mz:gsssg}). The only finite Coxeter groups $G$ with $\sigma ^0 (G) =0$ are 
the dihedral groups, $G _2$, $\Sigma _3$, $\Sigma _4$, and $D_3$. Also there are no finite Coxeter 
groups $G$ with $\sigma ^0 (G) =1$. 
In this paper, we will assume that $\sigma ^0 (G) >1$ for each group G that we 
are discussing. It is known that for groups with $\sigma ^0 (G) >1$, any 
generating pair will be a $(p,q,r)$ generating pair with $\frac{1}{p} +\frac{1}{q} + \frac{1}{r}<1$. 
Using the Riemann-Hurwitz equation, we see that given any generating pair of $G$, we get an upper bound 
on the strong symmetric genus of $G$ \cite{s:sfg}. If $G$ has a $(p,q,r)$ generating pair, then 
$\sigma ^0 (G)\leq 1+\frac{1}{2} |G|\cdot (1-\frac{1}{p} -\frac{1}{q} -\frac{1}{r})$. The following lemma, which is a result of 
Singerman \cite{s:srs} (see also \cite{mz:gsssg,t:fgas}), shows that the strong symmetric genus for many groups 
is computed directly from a minimal generating pair.

\begin{lemma}[Singerman \cite{s:srs}]\label{lemma:sing}
Let $G$ be a finite group such that $\sigma ^0 (G) >1$. If $|G|>12(\sigma ^0 (G) -1)$, then $G$ has a 
$(p,q,r)$ generating pair with 
\begin{equation}
\sigma ^0 (G)=1+\frac{1}{2} |G|\cdot (1-\frac{1}{p} -\frac{1}{q} -\frac{1}{r}).
\end{equation}
In addition, if $p\neq 2$, $p=q=3$, and $r$ is $4$ or $5$.
\end{lemma}

We also include a lemma that allows for control of minimal generating pairs when 
passing to quotient groups.

\begin{lemma}
Let $G$ be a finite group such that $G$ has a $(p,q,r)$ generating pair. For any normal subgroup $N \lhd G$,
any minimal generating pair of $G/N$ must have a $(p',q',r')$ generating pair such that 
$\frac{1}{p} +\frac{1}{q} + \frac{1}{r}\leq \frac{1}{p'} +\frac{1}{q'} + \frac{1}{r'}$. In other words, 
$\sigma ^0 (G/N)-1 \leq \frac{\sigma ^0 (G) - 1}{|N|} $.
\end{lemma}

Proof: Suppose that $(x,y)$ is a $(p,q,r)$ generating pair of $G$. In addition, let $\bar{x}$ and $\bar{y}$ be 
the images under the quotient map $\pi : G \rightarrow G/N$ of $x$ and $y$ respectively. Let $p'$, 
$q'$, and $r'$ be the orders of $\bar{x}$, $\bar{y}$, and $\bar{x}\bar{y}$. Clearly, $p'\leq p$, $q' \leq q$, 
and $r' \leq r$; therefore, $\frac{1}{p'} +\frac{1}{q'} + \frac{1}{r'}\geq \frac{1}{p} +\frac{1}{q} + \frac{1}{r}$. 
It is also clear that $\bar{x}$ and $\bar{y}$ generate $G/N$.
$\Box $

In order to use this quotient group result, we recall some cases where certain Coxeter groups are 
quotients of other Coxeter groups: For a fixed $n$, $\Sigma _n \cong B_n/(\mathbb{Z}_n)^n$, and 
$\Sigma _n \cong D_n/(\mathbb{Z}_n)^{n-1}$. In addition if $n$ is odd, $D_n \cong B_n / Z(B_n) $.

For an example of using quotient groups, we look at $D_{17}$. We notice that both $B_{17}$ and $\Sigma _{17}$ have 
minimal $(2,4,6)$ generating pairs. Now $D_{17}$ has a minimal generating pair; suppose $D_{17}$  
has a $(p,q,r)$ minimal generating pair. Since $D_{17} \cong B_{17}/Z(B_{17})$, $\frac{11}{12} \leq 
\frac{1}{p} +\frac{1}{q} + \frac{1}{r}$ with $p\leq 2$, $q\leq 4$, and $r\leq 6$. 
On the other hand, we have $\frac{1}{p} +\frac{1}{q} + \frac{1}{r}\leq \frac{11}{12}$, 
because $\Sigma _{17} \cong D_{17}/(\mathbb{Z}_2)^{16}$. So we see that $D_{17}$ must have a 
$(2,4,6)$ minimal generating pair.

\section{Generators of $D_n$ }\label{sec:gen}

For notation purposes, we recall that $D_n$ is an index two subgroup of $B_n$ and that $B_n$ is 
the wreath product $\mathbb{Z}_2 \wr \Sigma _n$. We will use notation that was adopted by V. S. Sikora 
\cite{s:tbhg} for the elements of $B_n$ and thus $D_n$  (see also \cite{j:ssgh} ). For an 
element of $B_n$, we will write a tuple $[\sigma , b ]$ where $\sigma $ is an element of $\Sigma _n$ and 
where $b$ is a list of $n$ binary digits representing the element of $(\mathbb{Z}_2)^n$. The multiplication then becomes 
$[\sigma , b ]\cdot [\tau  , c ]=[\sigma \cdot \tau , \tau ^{-1}(b)+c]$ where addition in the binary digits is a 
parity computation in each entry. We will use the convention of calling $b$ even or odd according to the number of 
ones appearing as binary digits of $b$. Notice that if $b$ and $c$ have the same parity, then $b+c$ is 
even; and if they differ in parity, then $b+c$ is odd. Using this notation, an element $[\sigma , b]\in B_n$ 
is contained in $D_n$ if and only if $b$ is even.

Notice that if $[\sigma , b]$ and $[\tau , c]$ generate $D_n$, then $\sigma $ and $\tau $ must generate $\Sigma _n$. 
Since we are looking to find generators of $D_n$, we need to first find generators of $\Sigma _n$. Next we 
construct generators of $D_n$ from the generators of $\Sigma _n$. Also 
during this construction we would like to have control of the orders of the generators of $D_n$ as well 
as to have control of the order of their product as described in Section \ref{sec:gp}. 

Recall that the following sequence is a split exact sequence of groups:
\begin{equation}
(\mathbb{Z}_2)^{n-1} \stackrel{i}{\rightarrow} D_n \stackrel{\pi}{\rightarrow} \Sigma _n 
\end{equation}
where $i(b)=[1,b]$ and $\pi ([\sigma ,b])=\sigma $. 
In addition, we have the following commutative diagram
\begin{equation}
\begin{CD} 
(\mathbb{Z}_2)^{n-1} @>{i }>> D_n @>{\pi }>> \Sigma _n \\
@VVV	                 @VVV		  @|\\
(\mathbb{Z}_2)^n @>{i }>> B_n @>{\pi }>> \Sigma _n \\
\end{CD}
\end{equation}
where both horizontal sequences are split exact. Using the bottom split exact sequence, the author has 
proven the following proposition:

\begin{proposition}[Jackson \cite{j:ssgh}]\label{prop:subbn}
For $n\geq 5$, let $G$ be a subgroup of $B_n$ with $\pi (G)=\Sigma _n$; then $G$ is a split extension of $\Sigma _n$ by one of 
the following: $1$, $\mathbb{Z}_2$, $(\mathbb{Z}_2)^{n-1}$, or $(\mathbb{Z}_2)^n$. In the first two cases, $G$ is 
isomorphic to $\Sigma _n$ and $\mathbb{Z}_2 \times \Sigma _n$, respectively;  in the third 
case either $$G=\{ [\sigma , b] \in B_n | b \textup{ is even} \}=D_n \textup{  or  } G=\left\{ [\sigma , b] \in B_n \left| \begin{array}{c}
b \textup{ is even if }\sigma \in A_n\\
b \textup{ is odd if }\sigma \in \Sigma _n \setminus A_n
\end{array} \right. \right\};$$ in the fourth case $G=B_n$.
\end{proposition}

Proposition \ref{prop:subbn} leads to Corollary \ref{cor:subdn}:

\begin{corollary}\label{cor:subdn}
For $n\geq 5$, let $H$ be a subgroup of $D_n$ with $\pi (H)=\Sigma _n$; then $H$ is a split extension of $\Sigma _n$ 
by one of the following: $1$, $\mathbb{Z}_2$, or $(\mathbb{Z}_2)^{n-1}$. In the first case, $H \cong \Sigma _n$. 
In the second case, which can only occur when $n$ is even,  $H \cong \Sigma _n \times Z(D_n)$. In the 
third case, $H=D_n$.
\end{corollary}

Using Corollary \ref{cor:subdn}, we can prove Proposition \ref{prop:gendn}, which we will use to construct 
generators of $D_n$ from those of $\Sigma _n$.

\begin{proposition}\label{prop:gendn}
Suppose $\sigma $ and $\tau $ generate $\Sigma _n$ as $Symm(\Gamma )$ where $\Gamma =\{ 1,2,\dots ,n\}$ 
with both $\sigma $ and $\sigma \cdot \tau $ having even order. Assume, furthermore, that $\sigma $ fixes two 
elements $i$ and $j$ of $\Gamma $, which are in the same cycle of the element $\sigma \cdot \tau $; assume as well 
that $\sigma $ fixes a third element of $\Gamma $ if $n$ is even. 
Let $b=(0,0,\dots , 0 ,0)$, and let $a=(0,\dots , 0,1,0,\dots ,0,1,0,\dots ,0)$ where there is a $1$ in the $i^{th}$ 
position and in the $j^{th}$ position. Under these conditions, $[\sigma , a]$ and 
$[\tau , b]$ generate $D_n$. In addition the elements $[\sigma , a]$, $[\tau , b]$ and $[\sigma \cdot \tau , \tau ^{-1} (a)]$
have the same orders as $\sigma $, $\tau $, and $\sigma \cdot \tau $, respectively.
\end{proposition}

Proof: The result concerning the orders of $[\sigma , a]$ and $[\tau , b]$ are obvious since $\sigma (a)=a$ and 
$\tau (b)=b$. On the other hand, $(\sigma \cdot \tau )^{-k}(\tau ^{-1}(a))=(\sigma \cdot \tau )^{-k-1}(a)$. 
Notice that if $k$ is the length of the cycle in $\sigma \cdot \tau $ that contains $i$ and $j$, then 
\begin{equation}
[\sigma \cdot \tau , \tau ^{-1}(a)]^k = [(\sigma \cdot \tau)^k , (0,0, \dots ,0,0)].
\end{equation}
It is clear then that the order of $[\sigma \cdot \tau , \tau ^{-1}(a) ]$ 
is the same as the order of $\sigma \cdot \tau $.

Let $H= \langle [\sigma , a], [\tau , b] \rangle \subset D_n$. We need to show that $H=D_n$. Notice that $H$ is a subgroup 
of $D_n$ such that $\pi (H)=\Sigma _n$. From Corollary \ref{cor:subdn} we know that $H$ is a split extension of 
$\Sigma _n$ by one of the following: $1$, $\mathbb{Z}_2$, or $(\mathbb{Z}_2)^{n-1}$; and we know that the second case only 
occurs when $n$ is even.

Recall that any section $s:\Sigma _n \rightarrow D_n$ takes $\alpha \in A_n$ to $[1,d]\cdot 
[\alpha ,(0,\dots ,0)]\cdot [1,d] ^{-1}$ for some $[1,d] \in B_n$; and if $\alpha \in \Sigma _n \setminus A_n$, 
then either $s(\alpha )=[1,d]\cdot [\alpha ,(0,\dots ,0)]\cdot [1,d] ^{-1}$ or  $s(\alpha )=[1,d]\cdot [\alpha ,(1,\dots ,1)]
\cdot [1,d] ^{-1}$. Also notice that the last case only occurs if $n$ is even.

We will show first that $[\sigma , a]$ cannot be in the image of any such section via contradiction.
Suppose $[\sigma , a]$ is in the image of some section homorphism $s:\Sigma _n \rightarrow D_n$. 
Now $[\sigma , a]$ cannot be $[1,d]\cdot [\sigma ,(0,\dots ,0)]\cdot [1,d] ^{-1}=[\sigma , \sigma ^{-1}(d)+d]$ 
for any $[1,d] \in B_n$ 
since $\sigma $ fixes $i$ and $j$ so that $\sigma ^{-1}(d)+d$ has a $0$ in both the $i^{th}$ and $j^{th}$ positions.
We, therefore, may assume that $n$ is even and 
\begin{equation}
[\sigma , a]=[1,d]\cdot [\sigma ,(1,\dots ,1)]\cdot [1,d] ^{-1}=
[\sigma ,  \sigma ^{-1}(d)+d + (1,\dots ,1) ]
\end{equation}
for some $[1,d] \in B_n$. $ \sigma ^{-1}(d)+d + (1,\dots ,1) \neq a $  since 
$\sigma $ fixes $k$ and $a$ has a zero in the $k^{th}$ position.

$[\sigma , a]$ is not in the image of any section; thus we only need to show that $[\sigma , a]$ is not 
equal to $s(\sigma ) \cdot [1,(1,\dots , 1)]$ for any section $s$. Notice that if $[\sigma, a]$ were such an element, 
then $a$ is either $\sigma ^{-1}(d)+d + (1,\dots ,1)$ or $\sigma ^{-1}(d)+d$, which are both 
ruled out in the previous paragraph. $\Box $

\section{Results}\label{sec:res}

From Section \ref{sec:gen}, we see that the strong symmetric genus of $D_n$ may be shown by demonstrating 
a particular generating pair of $\Sigma _n$. First we notice that if $\Sigma _n$ has a $(p,q,r)$ generating pair, 
then at most one of $p$, $q$, or $r$ is odd. This result also holds for generators of $D_n$. So we see that 
the best possible generating pair for any $D_n$ with $n>3$ is a $(2,3,8)$ generating pair. We notice 
that Marston Conder \cite{c:gasg} has demonstrated $(2,3,8)$ generating pairs for each $\Sigma _n$ with 
$n\geq 168$. If we call this generating pair $(\sigma , \tau )$, we notice that for each $n\geq 168$, excluding those for values 
of $n$ listed below, $\sigma $ fixes three elements, two of which are in the same cycle of $\sigma \cdot \tau $. 
In these cases we apply Proposition \ref{prop:gendn} to see that $D_n$ has a $(2,3,8)$ minimal generating pair.
The exceptional values of $n$ are 171, 173, 174, 181, 185, 188, 194, 201, 202, 206, 209, 214, 230, 250, 257, 265, 
and 286. 

In the remaining cases, the results were computed using GAP \cite{g:gap}. For $n\geq 30$, including the exceptional 
values of $n$ listed above, $D_n$ was shown to have a $(2,3,8)$ minimal generating pair. So for each $n\geq 30$, 
$\sigma ^0 (D_n)=\frac{n! \; 2^{n-1}}{48}$. For each $D_n$ with $n<30$ as well as for each sporadic group listed in Table I, 
an exhustive search was performed using GAP \cite{g:gap} 
to find a minimal generating pair. In this manner the new results found in Table I and Table II, as well as in 
Theorem \ref{thm:main}, were obtained.

\end{document}